\documentclass{psapm-l}
\usepackage{graphicx}

\newtheorem{theorem}{Theorem}[section]

\theoremstyle{definition}
\newtheorem{definition}[theorem]{Definition}

\theoremstyle{remark}
\newtheorem{remark}{Remark}[section]

\numberwithin{equation}{section}

\newcommand{\x}{{\bf x}}
\newcommand{\R}{\mathbb{R}}
\newcommand{\del}{{\partial}}
\newcommand{\vv}{\mathbf{v}}
\newcommand{\aaa}{\mathbf{a}}
\newcommand{\bb}{\mathbf{b}}
\newcommand{\xx}{{\bf x}}

\newcommand \alp{\alpha}
\newcommand \eps{\varepsilon}
\newcommand \vphi{\varphi}

\newcommand \gam{\gamma}

\newcommand \der{\partial}

\newcommand{\PtUpL}{{P_1}}
\newcommand{\PtLwL}{{P_2}}
\newcommand{\PtLwR}{{P_3}}
\newcommand{\PtUpR}{{P_4}}

\newcommand \Sonic{\Gamma_{sonic}}
\newcommand \Shock{\Gamma_{shock}}
\newcommand \Wedge{\Gamma_{wedge}}

\newcommand \inn{\text{in}}

\newcommand \Om{\Omega}

\begin{document}
\title[Shock Reflection-Diffraction]{Shock Reflection-Diffraction Phenomena\\
and\\ Multidimensional Conservation Laws}

\author{Gui-Qiang Chen}
\address{Department of Mathematics, Northwestern University, Evanston, IL 60208-2730, USA}
\email{gqchen@math.northwestern.edu}
\thanks{Gui-Qiang Chen's research was supported in part by the National
 Science Foundation under Grants DMS-0807551, DMS-0720925,
 and DMS-0505473.}
\author{Mikhail Feldman}
\address{Department of Mathematics, University of Wisconsin,
Madison, WI 53706-1388, USA}
\email{feldman@math.wisc.edu}
\thanks{Mikhail Feldman's research was
supported in part by the National Science Foundation under Grants
DMS-DMS-0800245
and DMS-0354729.}

\date{October 15, 2008}


\begin{abstract}
When a plane shock hits a wedge head on, it experiences a
reflection-diffraction process, and then a self-similar reflected
shock moves outward as the original shock moves forward in time. The
complexity of reflection-diffraction configurations was first
reported by Ernst Mach in 1878, and experimental, computational, and
asymptotic analysis has shown that various patterns of shock
reflection-diffraction configurations may occur, including regular
reflection and Mach reflection. In this paper we start with various
shock reflection-diffraction phenomena, their fundamental scientific
issues, and their theoretical roles as building blocks and
asymptotic attractors of general solutions in the mathematical
theory of multidimensional hyperbolic systems of conservation laws.
Then we describe how the global problem of shock
reflection-diffraction by a wedge can be formulated as a free
boundary problem for nonlinear conservation laws of mixed-composite
hyperbolic-elliptic type. Finally we discuss some recent
developments in attacking the shock reflection-diffraction problem,
including the existence, stability, and regularity of global regular
reflection-diffraction solutions. The approach includes techniques
to handle free boundary problems, degenerate elliptic equations, and
corner singularities, which is highly motivated by experimental,
computational, and asymptotic results. Further trends and open
problems in this direction are also addressed.
\end{abstract}
\maketitle

\section{Introduction}

Shock waves occur in many physical situations in nature. For
example, shock waves can be produced by solar winds (bow shocks),
supersonic or near sonic aircrafts (transonic shocks around the
body), explosions (blast waves), and various natural processes. When
such a shock hits an obstacle (steady or flying), shock
reflection-diffraction phenomena occur. One of the most important
problems in mathematical fluid dynamics is the problem of shock
reflection-diffraction by a wedge. When the plane shock hits a wedge
head on, it experiences a reflection-diffraction process, and then a
fundamental question is what types of wave patterns of
reflection-diffraction configurations it may form around the wedge.

The complexity of reflection-diffraction configurations was first
reported by Ernst Mach \cite{Mach} in 1878, who first observed two
patterns of reflection-diffraction configurations: regular
reflection (two-shock configuration) and Mach reflection
(three-shock configuration); also see
\cite{BD,Krehl-Greest,Reichenback}. The problem remained dormant
until the 1940's when von Neumann, Friedrichs, Bethe, as well as
many experimental scientists, among others, began extensive research
into all aspects of shock reflection-diffraction phenomena, due to
its importance in applications. See von Neumann
\cite{Neumann1,Neumann2} and Ben-Dor \cite{BD};  also see
\cite{Bazhennova-Gvozdeva,Glass,HOS,IVF,Korobeinikov,Lyakhov-Podlubny,Srivastava1,Srivastava2}
and the references cited therein.
It has been found that the situation is much more complicated than
what Mach originally observed: The Mach reflection can be further
divided into more specific sub-patterns, and various other patterns
of shock reflection-diffraction may occur such as the von Neumann
reflection and the Guderley reflection; see
\cite{BD,CF,GlimmMajda,Guderley,HT,SA,TH,TSK1,TSK2,VD,Neumann1,Neumann2}
and the references cited therein.

Then the fundamental scientific issues include:

(i) {\it Structure of the shock reflection-diffraction
configurations};

(ii) {\it Transition criteria between the different patterns of
shock reflection-diffraction configurations};

(iii) {\it Dependence of the patterns upon the physical parameters
such as the wedge angle $\theta_w$, the incident-shock-wave Mach
number $M_s$, and the adiabatic exponent $\gamma\ge 1$}.

Careful asymptotic analysis has been made for various
reflection-diffraction configurations in  Lighthill
\cite{Lighthill}, Keller-Blank \cite{KB}, Hunter-Keller \cite{HK},
Morawetz \cite{Morawetz2},
\cite{GRT,TR,Harabetian,hunter1,Neumann1,Neumann2}, and the
references cited therein; also see Glimm-Majda \cite{GlimmMajda}.
Large or small scale numerical simulations have been also made; see,
e.g. \cite{BD,GlimmMajda},
\cite{DP,DG,Hindman-Kutler-Anderson,Kutler-Shankar,Schneyer,Shankar-Kutler-Anderson},
\cite{Berger-Colella,Colella-Glaz,Glaz-Colella1,Glaz-Colella2,GWGH,IVF,WC},
and the references cited therein.

However, most of the fundamental issues for shock
reflection-diffraction phenomena have not been understood,
especially the global structure and transition of different patterns
of shock reflection-diffraction configurations. This is partially
because physical and numerical experiments are hampered by various
difficulties and have not been able to select the correct transition
criteria between different patterns. In particular, numerical
dissipation or physical viscosity smear the shocks and cause
boundary layers that interact with the reflection-diffraction
patterns and can cause {\it spurious Mach steams}; cf.
Woodward-Colella \cite{WC}. Furthermore, some difference between two
different patterns are only fractions of a degree apart (e.g., see
Fig. 5 below), a resolution even by sophisticated modern experiments
(e.g. \cite{LD}) has been unable to reach, as pointed out by Ben-Dor
in \cite{BD}: ``{\it For this reason it is almost impossible to
distinguish experimentally between the sonic and detachment
criteria}" (cf. Section 5 below). In this regard, it seems that the
ideal approach to understand fully the shock reflection-diffraction
phenomena, especially the transition criteria, is still via rigorous
mathematical analysis. To achieve this, it is essential to establish
first the global existence, regularity, and structural stability of
solutions of the shock reflection-diffraction problem.

On the other hand,  shock reflection-diffraction configurations are
the core configurations in the structure of global solutions of the
two-dimensional Riemann problem for hyperbolic conservation laws;
while the Riemann solutions are building blocks and local structure
of general solutions and determine global attractors and asymptotic
states of entropy solutions, as time tends to infinity, for
multidimensional hyperbolic systems of conservation laws. See
\cite{ChangChen,CCY1,CCY2,CH,GlimmMajda,GlimmK,KTa,LaxLiu,LZY,SCG,Ser3,Serre,Zhe}
and the references cited therein. In this sense, we have to
understand the shock reflection-diffraction phenomena, in order to
understand fully entropy solutions to multidimensional hyperbolic
systems of conservation laws.

In this paper, we first formulate the shock reflection-diffraction
problem into an initial-boundary value problem in Section 2. Then we
employ the essential feature of self-similarity of the
initial-boundary value problem to reformulate the problem into a
boundary value problem in the unbounded domain in Section 3. In
Section 4, we present the unique solution of normal reflection for
this problem when the wedge angle is ${\pi}/{2}$. In Section 5, we
exhibit the local theory of regular reflection-diffraction,
introduce a stability criterion to determine state (2) at the
reflection point on the wedge, and present the von Neumann's
detachment and sonic conjectures. Then we discuss the role of the
potential flow equation in the shock reflection-diffraction problem
even in the level of the full Euler equations in Section 6. Based on
the local theory, we reduce the boundary value problem into a free
boundary problem in the context of potential flow in Section 7. In
Section 8, we describe a global theory for regular
reflection-diffraction for potential flow, established in
Chen-Feldman \cite{ChenFeldman2a,ChenFeldman2b,ChenFeldman3} and
Bae-Chen-Feldman \cite{BCF}. In Section 9, we discuss some open
problems and new mathematics required for further developments,
which are also essential for solving multidimensional problems in
conservation laws and other areas in nonlinear partial differential
equations.

\section{Mathematical Formulation I: Initial-Boundary Value Problem}

The full Euler equations for compressible fluids in
$\R^{3}_+:=\R_+\times \R^2, t\in\R_+:=(0,\infty), \x\in\R^2$, are of
the following form:
\begin{equation}\label{E-1}
\left\{\begin{aligned} &\del_t\,\rho +\nabla_\x\cdot (\rho\vv)=0,\\
&\del_t(\rho\vv)+\nabla_\x\cdot\left(\rho \vv\otimes\vv\right)
  +\nabla p=0,\\
&\del_t(\frac{1}{2}\rho |\vv|^2+\rho
e)+\nabla_\x\cdot\big((\frac{1}{2}\rho|\vv|^2+\rho e + p)\vv\big)=0,
\end{aligned}
\right.
\end{equation}
where $\rho$ is the density, $\vv=(u, v)$ the fluid velocity, $p$
the pressure, and $e$ the internal energy. Two other important
thermodynamic variables are  the temperature $\theta$ and the energy
$S$.  The notation $\aaa\otimes\bb$ denotes the tensor product of
the vectors $\aaa$ and $\bb$.

Choosing $(\rho, S)$ as the independent thermodynamical variables,
then the constitutive relations can be written as
$(e,p,\theta)=(e(\rho,S), p(\rho,S),\theta(\rho,S))$ governed by
$$
\theta dS=de +pd\tau=de-\frac{p}{\rho^2}d\rho.
$$

For a polytropic gas, \begin{equation}\label{gas-1} p=(\gamma-1)\rho
e, \qquad e=c_v\theta, \qquad \gamma=1+\frac{R}{c_v},
\end{equation}
or equivalently,
\begin{equation}\label{gas-2}
p=p(\rho,S)=\kappa\rho^\gamma e^{S/c_v}, \qquad e=e(\rho,
S)=\frac{\kappa}{\gamma-1}\rho^{\gamma-1}e^{S/c_v},
\end{equation}
where $R>0$ may be taken to be the universal gas
    constant divided by the effective molecular weight of the particular
    gas,
    $c_v>0$ is the specific heat at constant volume,
    $\gamma>1$ is the adiabatic exponent, and $\kappa>0$ is any constant under
    scaling.

When a flow is potential, that is, there is a velocity potential
$\Phi$ such that
$$
\vv=\nabla_\x\Phi,
$$
then the Euler equations for the flow take the form:
\begin{equation}\label{Euler7}
\left\{\begin{aligned}
&\del_t \rho +\text{div}(\rho\nabla_\x\Phi)=0,\qquad\qquad\quad \text{(conservation of mass)}\\
&\del_t\Phi +\frac{1}{2}|\nabla_\x\Phi|^2 + i(\rho)=B_0, \qquad
\text{(Bernoulli's law)}
\end{aligned}\right.
\end{equation}
where
$$
i(\rho)=\frac{\rho^{\gamma-1}-1}{\gamma-1} \qquad\text{when}\,
\gamma>1,
$$
especially, $i(\rho)=\ln\rho$ when $\gamma=1$, by scaling and $B_0$
is the Bernoulli constant, which is usually determined by the
boundary conditions if such conditions are prescribed. {}From the
second equation in \eqref{Euler7}, we have
\begin{equation}\label{2.5-a}
\rho(D\Phi)=i^{-1}\big(B_0-(\del_t\Phi
+\frac{1}{2}|\nabla_\x\Phi|^2)\big).
\end{equation}
Then system
\eqref{Euler7} can be rewritten as the following time-dependent
potential flow equation of second order:
\begin{equation}\label{Euler8}
\del_t\rho(D\Phi) +\nabla\cdot (\rho(D\Phi)\nabla \Phi)=0
\end{equation}
with \eqref{2.5-a}.

\medskip
For a steady solution $\Phi=\varphi(\xx)$, i.e., $\del_t\Phi=0$, we
obtain the celebrated steady potential flow equation in
aerodynamics:
\begin{equation}\label{Euler9}
\nabla_\x\cdot (\rho(\nabla_\x\Phi)\nabla_\x\Phi)=0.
\end{equation}

\medskip
In applications in aerodynamics, \eqref{Euler7} or \eqref{Euler8} is
used for discontinuous solutions, and the empirical evidence is that
entropy solutions of \eqref{Euler7} or \eqref{Euler8} are fairly
close to entropy solutions for \eqref{E-1}, provided the shock
strengths are small, the curvature of shock fronts is not too large,
and the amount of vorticity is small in the region of interest.
Furthermore, we will show in Section 6 that, for the shock
reflection-diffraction problem, the Euler equations for potential
flow is actually {\it exact} in an important region of the solution
(see Theorem 6.1 below).

Then the problem of shock reflection-diffraction by a wedge can be
formulated as follows:

\medskip
{\bf Problem 2.1 (Initial-boundary value problem)}. {\it Seek a
solution of system \eqref{E-1} satisfying the initial condition at
$t=0$:
\begin{equation} \label{ibv-1}
(\vv, p,\rho) =\begin{cases}
(0,0,p_0,\rho_0), &\, \quad |x_2|>x_1\tan\theta_w, x_1>0,\\
(u_1,0,p_1,\rho_1), &\, \quad x_1<0;
\end{cases}
\end{equation}
and the slip boundary condition along the wedge boundary:
\begin{equation}\label{boundary-condition}
\vv\cdot \nu=0,
\end{equation}
where ${\bf \nu}$ is the exterior unit normal to the wedge boundary,
and state $(0)$ and $(1)$ satisfy
\begin{equation}\label{2.10}
u_1=\sqrt{\frac{(p_1-p_0)(\rho_1-\rho_0)}{\rho_0\rho_1}}, \qquad
\frac{p_1}{p_0}=\frac{(\gamma+1)\rho_1-(\gamma-1)\rho_0}{(\gamma+1)\rho_0-(\gamma-1)\rho_1},
\qquad \rho_1>\rho_0.
\end{equation}
That is, given $\rho_0, p_0, \rho_1$, and $\gamma>1$, the other
variables $u_1$ and $p_1$ are determined by \eqref{2.10}. In
particular, the Mach number $M_1=u_1/c_1$ is determined by
\begin{equation}\label{2.11}
M_1^2=\frac{2(\rho_1-\rho_0)^2}{\rho_0\big((\gamma+1)\rho_1-(\gamma-1)\rho_0\big)},
\end{equation}
where $c_1=\sqrt{\gamma p_1/\rho_1}$ is the sonic speed of fluid
state $(1)$. }

 \begin{figure}[h]
 \centering
 \includegraphics[height=2.0in,width=1.8in]{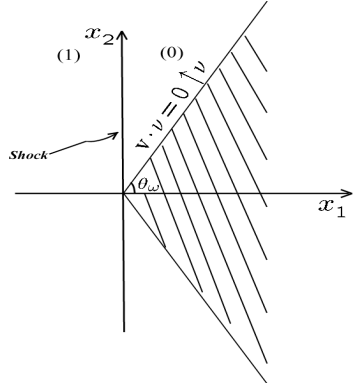}  
 \caption[]{Initial-boundary value problem}\label{fig:IBV-1}
 \end{figure}

\section{Mathematical Formulation II: Boundary Value Problem}

Notice that the initial-boundary value problem ({\bf Problem 2.1})
is invariant under the self-similar scaling:
$$
(t, \x)\longrightarrow (\alpha t, \alpha \x) \qquad\mbox{for any}\,
\, \alpha\ne 0.
$$
Therefore, we seek self-similar solutions:
$$
(\vv, p,\rho)(t,\x)=(\vv,p,\rho)(\xi,\eta), \qquad
(\xi,\eta)=\frac{\x}{t}.
$$
Then the self-similar solutions are governed by the following
system:
\begin{equation}\label{Euler-sms-1}
\left\{\begin{aligned} &(\rho U)_\xi +(\rho V)_\eta +2\rho=0,\\
&(\rho U^2+p)_\xi +(\rho UV)_{\eta}+3\rho
U=0,\\
&(\rho UV)_\xi +(\rho V^2+p)_{\eta}+3\rho
V=0,\\
&\big(U(\frac{1}{2}\rho q^2+\frac{\gamma p}{\gamma-1})\big)_\xi
+\big(V(\frac{1}{2}\rho q^2+\frac{\gamma p}{\gamma-1})\big)_\eta
+2(\frac{1}{2}\rho q^2+\frac{\gamma p}{\gamma-1})=0,
\end{aligned}
\right.
\end{equation}
where $q=\sqrt{U^2+V^2}$, and $(U,V)=(u-\xi, v-\eta)$ is the
pseudo-velocity.

The eigenvalues of system \eqref{Euler-sms-1} are
$$
\lambda_0=\frac{V}{U}\,\, \mbox{(repeated)},\qquad\,
\lambda_\pm=\frac{UV\pm c\sqrt{q^2-c^2}}{U^2-c^2},
$$
where $c=\sqrt{\gamma p/\rho}$ is the sonic speed.

When the flow is pseudo-subsonic, i.e.,  $q<c$,  the eigenvalues
$\lambda_\pm$ become complex and thus the system consists of two
transport equations and two nonlinear equations of
hyperbolic-elliptic mixed type. Therefore, system
\eqref{Euler-sms-1} is {\em hyperbolic-elliptic composite-mixed} in
general.

Since the problem is symmetric with respect to the axis $\eta=0$, it
suffices to consider the problem in the half-plane $\eta>0$ outside
the half-wedge:
$$
\Lambda:=\{\xi<0,\eta>0\}\cup\{\eta>\xi \tan\theta_w,\, \xi>0\}.
$$
Then the initial-boundary value problem ({\bf Problem 2.1}) in the
$(t, {\bf x})$--coordinates can be formulated as the following
boundary value problem in the self-similar coordinates $(\xi,\eta)$:

\medskip
{\bf Problem 3.1 (Boundary value problem in the unbounded domain)}.
{\it Seek a solution to system \eqref{Euler-sms-1} satisfying the
slip boundary condition on the wedge boundary and the matching
condition on the symmetry line $\eta=0$:
\begin{equation*}
(U,V)\cdot\nu =0 \qquad\,\, \hbox{on } \,\,
\partial\Lambda=\{\xi\le 0, \eta=0\}\cup\{\xi>0, \eta\ge
\xi\tan\theta_w\},
\end{equation*}
the asymptotic boundary condition as $\xi^2+\eta^2\to \infty$:
\begin{equation*}
{(U+\xi,V+\eta,p,\rho)\longrightarrow {\begin{cases}(0,
0,p_0,\rho_0),\qquad
                         \xi>\xi_0, \eta>\xi\tan\theta_w,\\
              (u_1, 0, p_1,\rho_1),\qquad
                          \xi<\xi_0, \eta>0.
\end{cases}}
}
\end{equation*}
}

 \begin{figure}[h]
 \centering
 \includegraphics[height=2.1in,width=2.3in]{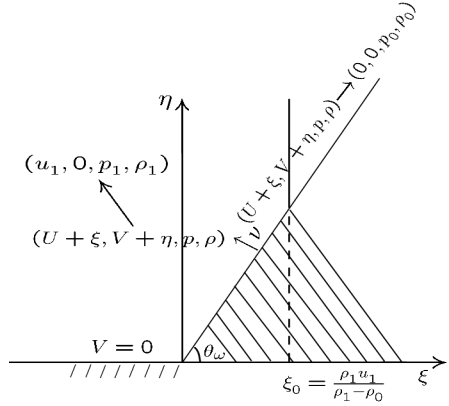}  
 \caption[]{Boundary value problem in the unbounded domain $\Lambda$}\label{fig:PM-1}
 \end{figure}

It is expected that the solutions of {\bf Problem 3.1} contain all
possible patterns of shock reflection-diffraction configurations as
observed in physical and numerical experiments; cf.
\cite{BD,CF,GlimmMajda,Guderley,Korobeinikov,Lyakhov-Podlubny,Mach,SA,TSK2}
and the references cited therein.

\section{Normal Reflection}

\begin{figure}[h]
 \centering
 \includegraphics[height=1.9in,width=2.9in]{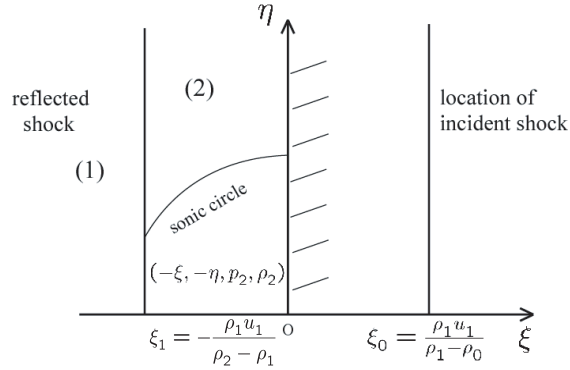}  
 \caption[]{Normal reflection solution}\label{fig:NR}
 \end{figure}

The simplest case of the shock reflection-diffraction problem is
when the wedge angle $\theta_w$ is $\pi/2$. In this case, the
reflection-diffraction problem simply becomes the normal reflection
problem, for which the incident shock normally reflects, and the
reflected shock is also a plane. It can be shown that there exist a
{\em unique} state $(p_2,\rho_2), \rho_2>\rho_1,$ and a {\em unique}
location of the reflected shock
\begin{equation}\label{4.1}
 \xi_1=-\frac{\rho_1 u_1}{\rho_2-\rho_1}\quad \qquad\text{with}\quad
 u_1=\sqrt{\frac{(p_2-p_1)(\rho_2-\rho_1)}{\rho_1\rho_2}}
\end{equation}
such that state $(2)=(-\xi, -\eta, p_2,\rho_2)$ is subsonic inside
the sonic circle with center at the origin and radius
$c_2=\sqrt{\gamma p_2/\rho_2}$, and is supersonic outside the sonic
circle (see Fig. \ref{fig:NR}). That is, in this case, the normal
reflection solution is unique.

\bigskip
In this case,
\begin{equation}\label{4.2}
 M_1^2=
\frac{2(\rho_2-\rho_1)^2}{\rho_2\big((\gamma+1)\rho_1-(\gamma-1)\rho_2\big)},
\qquad
\frac{p_2}{p_1}=\frac{(\gamma+1)\rho_2-(\gamma-1)\rho_1}{(\gamma+1)\rho_1-(\gamma-1)\rho_2},
\end{equation}
and $\frac{\rho_2}{\rho_1}=t>1$ is the unique root of
$$
\big(1+\frac{\gamma-1}{2}M_1^2\big)t^2-\big(2+\frac{\gamma+1}{2}M_1^2\big)t
+1=0,
$$
that is,
\begin{equation}\label{4.3}
 \frac{\rho_2}{\rho_1}
=\frac{4+(\gamma+1)M_1^2+M_1\sqrt{16+(\gamma+1)^2M_1^2}}{2\big(2+(\gamma-1)M_1^2\big)}.
\end{equation}
In other words, given $\rho_0, p_0, \rho_1$, and $\gamma>1$, state
$(2)=(-\xi, -\eta, p_2,\rho_2)$ is uniquely determined through
\eqref{2.10}--\eqref{2.11} and \eqref{4.1}--\eqref{4.3}.

\section{Local Theory and von Neumann's Conjectures
for Regular Reflection-Diffraction Configuration}

For a wedge angle $\theta_w\in (0, \pi/2)$, different
reflection-diffraction patterns may occur. Various criteria and
conjectures have been proposed for the existence of configurations
for the patterns (cf. Ben-Dor \cite{BD}). One of the most important
conjectures made by von Neumann \cite{Neumann1,Neumann2} in 1943 is
the {\em detachment conjecture}, which states that the regular
reflection-diffraction configuration may exist globally whenever the
two shock configuration (one is the incident shock and the other the
reflected shock) exists locally around the point $P_0$ (see Fig. 4).

The following theorem was rigorously shown in Chang-Chen
\cite{ChangChen} (also see Sheng-Yin \cite{ShengYin}, Bleakney-Taub
\cite{BT}, Neumann \cite{Neumann1,Neumann2}).

\begin{theorem}[Local theory]\label{local}
There exists $\theta_d=\theta_d(M_s,\gamma)\in (0, \pi/2)$ such
that, when $\theta_w\in (\theta_d,\pi/2)$, there are two states
$(2)=(U_2^a,V_2^a,p_2^a,\rho_2^a)$ and
$(U_2^b,V_2^b,p_2^b,\rho_2^b)$ such that
$$
|(U_2^a, V_2^a)|>|(U_2^b, V_2^b)| \quad \text{and}\quad |(U_2^b,
V_2^b)| < c_2^b,
$$
where $c_2^b=\sqrt{\gamma p^b_2/\rho^b_2}$ is the sonic speed.
\end{theorem}

Then the conjecture can be stated as follows:

\medskip
{\bf The von Neumann's Detachment Conjecture}: {\em There exists a
global regular reflection-diffraction configuration whenever the
wedge angle $\theta_w$ is in $(\theta_d, \pi/2)$}.

\medskip
It is clear that the regular reflection-diffraction configuration is
not possible without a local two-shock configuration at the
reflection point on the wedge, so this is the weakest possible
criterion. In this case, the local theory indicates that there are
two possible states for state (2). There had been a long debate to
determine which one is more physical for the local theory; see
Courant-Friedrichs \cite{CF}, Ben-Dor \cite{BD}, and the references
cited therein.

Since the reflection-diffraction problem is not a local problem, we
take a different point of view that the selection of state (2)
should be determined by the global features of the problem, more
precisely, by the stability of the configuration with respect to the
wedge angle $\theta_w$, rather than the local features of the
problem.

\medskip
{\bf Stability Criterion to Select the Correct State (2)}: {\it
Since the solution is unique when the wedge angle $\theta_w=\pi/2$,
it is required that our global regular reflection-diffraction
configuration should be stable and converge to the unique normal
reflection solution when $\theta_w\to \pi/2$, provided that such a
global configuration can be constructed.}

\medskip
We employ this stability criterion to conclude that our choice for
state (2) must be $(U_2^a, V_2^a, p_2^a, \rho_2^a)$. In general,
$(U_2^a, V_2^a, p_2^a, \rho_2^a)$ may be supersonic or subsonic. If
it is supersonic, the propagation speeds are finite and state (2) is
completely determined by the local information: state (1), state
(0), and the location of the point $P_0$. This is, any information
from the reflected region, especially the disturbance at the corner
$P_3$, cannot travel towards the reflection point $P_0$. However, if
it is subsonic, the information can reach $P_0$ and interact with
it, potentially altering the reflection-diffraction type. This
argument motivated the second conjecture as follows:

\medskip
{\bf The von Neumann's Sonic Conjecture}: {\em There exists a
regular reflection-diffraction configuration when $\theta_w\in
(\theta_s, \pi/2)$ for $\theta_s>\theta_d$ such that $|(U_2^a,
V_2^a)|>c_2^a$ at $P_0$.}
\begin{figure}[h]
 \centering
 \includegraphics[height=1.8in,width=2.2in]{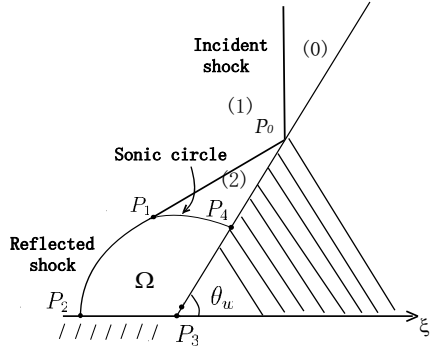}  
 \caption[]{Regular reflection-diffraction configuration}\label{fig:RC-1}
 \end{figure}

This sonic conjecture is based on the following fact: If state (2)
is sonic when $\theta_w=\theta_s$, then $|(U^a_2, V^a_2)|>c_2^a$ for
any $\theta_w\in (\theta_s, \pi/2)$. This conjecture is stronger
than the detachment one. In fact, the regime between the angles
$\theta_s$ and $\theta_d$ is very narrow and is only fractions of a
degree apart; see Fig. 5 from Sheng-Yin \cite{ShengYin}.

\begin{figure}[h]
 \centering
 \includegraphics[height=2.0in,width=2.2in]{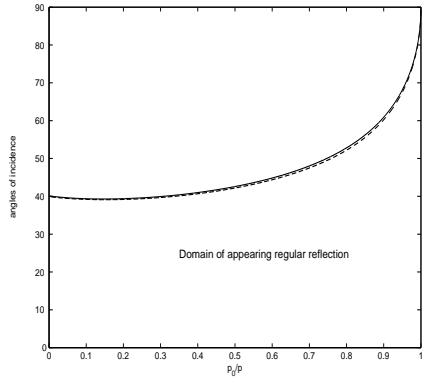}  
 \caption[]{The von Neumann's sonic criterion vs the detachment criterion $\theta_s>\theta_d$
 when $\gamma=1.4$.}
 \label{fig:DS-1}
 \end{figure}

\section{The Potential Flow Equation}

In this section, we discuss the role of the potential flow equation
in the shock reflection-diffraction problem for the full Euler
equations.

Under the Hodge-Helmoltz decomposition $(U,V)=\nabla\varphi +W$ with
$\nabla\cdot W=0$, the Euler equations \eqref{Euler-sms-1} become
\begin{eqnarray}
&&\nabla\cdot (\rho\nabla\varphi)+2\rho + \nabla\cdot(\rho W)=0,\label{6.1}\\
&&\nabla(\frac{1}{2}|\nabla\varphi|^2+\varphi)+\frac{1}{\rho}\nabla
p =(\nabla\varphi +W)\cdot\nabla W+
(\nabla^2\varphi+I)W,\label{6.2}\\
&&(\nabla\varphi+W)\cdot \nabla \omega +(1+\Delta
\varphi)\omega=0,\label{6.3}\\
&&(\nabla\varphi+W)\cdot \nabla S=0,\label{6.4}
\end{eqnarray}
where $\omega={\rm curl}\, W={\rm curl} (U,V)$ is the vorticity of
the fluid, $S=c_v\, {\rm ln} (p\rho^{-\gamma})$ is the entropy, and
the gradient $\nabla$ is with respect to the self-similar variables
$(\xi, \eta)$ from now on.

When $\omega=0, S=const.$, and $W=0$ on a curve $\Gamma$ transverse
to the fluid direction, we first conclude from \eqref{6.3} that, in
the domain $\Omega_1$ determined by the fluid trajectories:
$$
\frac{d}{dt}(\xi,\eta)=(\nabla\varphi+W)(\xi,\eta)
$$
past $\Gamma$,
$$
\omega=0, \qquad \text{i.e.}\quad \text{curl}\, W=0.
$$
This implies that $W=const.$ since $\nabla\cdot W=0$. Then we
conclude that
$$
W=0   \qquad \text{in}\,\, \Omega_1,
$$
since $W|_\Gamma=0$, which yields that the right-hand side of
equation \eqref{6.2} vanishes. Furthermore, from \eqref{6.4},
$$
S=const. \qquad\text{in}\,\, \Omega_1,
$$
which implies that
$$
p=const.\, \rho^\gamma.
$$
By scaling, we finally conclude that the solutions of system
\eqref{6.1}--\eqref{6.4} in the domain $\Omega_1$ is determined by
the following system for self-similar solutions:
\begin{equation}\label{potential-eq}
\left\{\begin{array}{ll} \nabla\cdot (\rho\nabla\varphi)+2\rho=0,\\
\frac{1}{2}|\nabla\varphi|^2+\varphi
+\frac{\rho^{\gamma-1}}{\gamma-1}=\frac{\rho_0^{\gamma-1}}{\gamma-1}.
\end{array}\right.
\end{equation}
or the potential flow equation for self-similar solutions:
\begin{equation}\label{potential-1}
\nabla\cdot \big(\rho(\nabla\varphi, \varphi)\nabla\varphi\big)
+2\rho(\nabla\varphi, \varphi)=0,
\end{equation}
with
\begin{equation}
\rho(|\nabla\varphi|^2, \varphi)
=\big(\rho_0^{\gamma-1}-(\gamma-1)(\varphi+\frac{1}{2}|\nabla\varphi|^2)\big)^{\frac
1{\gamma-1}}. \label{1.1.6}
\end{equation}
Then we have
\begin{equation}\label{c-through-density-function}
c^2=c^2(|\nabla\varphi|^2,\varphi,\rho_0^{\gamma-1})
=\rho_0^{\gamma-1}-(\gamma-1)(\frac{1}{2}|\nabla\varphi|^2+\varphi).
\end{equation}

For our problem (see Fig. 4), we note that, for state (2),
\begin{equation}\label{across-sonic}
\omega=0, \quad W=0, \quad S=S_2.
\end{equation}
Then,  if our solution $(U,V,p,\rho)$ is $C^{0,1}$  and the gradient
of the tangential component of the velocity is continuous across the
sonic arc $\Gamma_{sonic}$,
we still have \eqref{across-sonic} along $\Gamma_{sonic}$ on the
side of $\Omega$. Thus, we have

\begin{theorem}\label{p-dominate}
Let $(U,V,p,\rho)$ be a solution of our {\bf Problem 3.1} such that
$(U,V,p,\rho)$  is $C^{0,1}$ in the open region $P_0P_1P_2P_3$ and
the gradient of the tangential component of $(U,V)$ is continuous
across the sonic arc $\Sonic$.
Let $\Omega_1$ be the subregion of $\Omega$ formed by the fluid
trajectories past the sonic arc $\Gamma_{sonic}$, then, in
$\Omega_1$, the potential flow equation \eqref{potential-1} with
\eqref{1.1.6} coincides with the full Euler equations
\eqref{6.1}--\eqref{6.4}, that is, equation \eqref{potential-1} with
\eqref{1.1.6} is exact in the domain $\Omega_1$ for {\bf Problem
3.1}.
\end{theorem}

\begin{remark}
The regions such as $\Omega_1$ also exist in various Mach
reflection-diffraction configurations. Theorem \ref{p-dominate}
applies to such regions whenever the solution $(U,V,p,\rho)$ is
$C^{0,1}$ and the gradient of the tangential component of $(U,V)$ is
continuous. In fact, Theorem 8.3 indicates that, for the solutions
$\varphi$ of \eqref{potential-1} with \eqref{1.1.6}, the $C^{1,1}$
regularity of $\varphi$ and the continuity of the tangential
component of the velocity field $(U,V)=\nabla\varphi$ are optimal
 across the sonic arc $\Gamma_{sonic}$.
\end{remark}

\begin{remark}
The importance of the potential flow equation \eqref{Euler8} with
\eqref{2.5-a} in the time-dependent Euler flows was also observed by
Hadamard \cite{Hadamard} through a different argument.
\end{remark}

Furthermore, when the wedge angle $\theta_w$ is close to $\pi/2$, it
is expected that the curvature of the reflected shock is small so
that, in the other part $\Omega_2$ of $\Omega$, the vorticity
$\omega$ is small and the entropy is close to the constant. Then, in
the reflection-diffraction domain $\Omega=\Omega_1\cup \Omega_2$,
the potential flow equation \eqref{potential-1} with \eqref{1.1.6}
dominates, provided that the exact state along the reflected shock
is given.

\begin{figure}[h]
 \centering
\includegraphics[height=1.7in,width=2.4in]{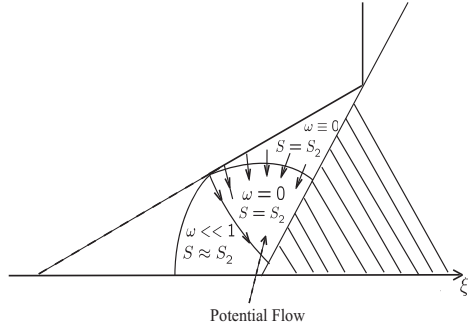}  
\caption[]{The potential flow equation dominates the domain
$\Omega$}
 \label{fig:PF-dominate}
 \end{figure}
\medskip

Equation \eqref{potential-1} with \eqref{1.1.6} is a nonlinear
equation of mixed elliptic-hyperbolic type. It is elliptic if and
only if
\begin{equation}
|\nabla\varphi| < c(|\nabla\varphi|^2,\varphi,\rho_0^{\gamma-1}),
\label{1.1.8}
\end{equation}
which is equivalent to
\begin{equation}
|\nabla \varphi| <c_*(\varphi, \rho_0, \gamma)
:=\sqrt{\frac{2}{\gamma+1}\big(\rho_0^{\gamma-1}-(\gamma-1)\varphi\big)}.
\label{1.1.8a}
\end{equation}

\bigskip
The study of partial differential equations of mixed
hyperbolic-elliptic type can date back 1940s (cf.
\cite{Bet,Chen0,Chen1,Wu,Yang}). Linear models of partial
differential equations of mixed hyperbolic-elliptic type include the
Lavrentyev-Betsadze equation:
$$
\partial_{xx}u+\text{sign}(x)\partial_{yy}u=0,
$$
the Tricomi equation:
$$
u_{xx}+ x u_{yy}=0 \qquad  (\text{hyperbolic degeneracy at} \,\,
x=0),
$$
and  the Keldysh equation:
$$
x u_{xx}+u_{yy}=0 \qquad (\text{parabolic degeneracy at}\,\, x=0).
$$
Nonlinear models of mixed-type equations for \eqref{potential-1}
with \eqref{1.1.6} include the {transonic small disturbance
equation}:
$$
\big((u-x)u_x+\frac{u}{2}\big)_x +u_{yy}=0
$$
or, for $v=u-x$,
$$
(v\, v_{x})_x+v_{yy} + \frac{3}{2}v_x+\frac{1}{2}=0,
$$
which has been studied in
\cite{CKK1,CanicKeyfitz,hunter1,hunter2,hunter3,Morawetz2} and the
references cited therein. Also see \cite{CKK2,Song,Zhe1,Zheng1} for
the models for self-similar solutions from the pressure gradient
system and nonlinear wave equations.

\section{Mathematical Formulation III: Free Boundary Problem for
Potential Flow}

For the potential equation \eqref{potential-1} with \eqref{1.1.6},
shocks are discontinuities in the pseudo-velocity $\nabla\varphi$.
That is, if $D^+$ and $D^-:=D\setminus\overline{D^+}$ are two
nonempty open subsets of $D\subset\R^2$ and $S:=\partial D^+\cap D$
is a $C^1$--curve where $D\varphi$ has a jump, then $\varphi\in
W^{1,1}_{loc}(D)\cap C^1(D^\pm\cup S)\cap C^2(D^\pm)$ is a global
weak solution of (\ref{potential-1}) with \eqref{1.1.6} in $D$ if
and only if $\varphi$ is in $W^{1,\infty}_{loc}(D)$ and satisfies
equation \eqref{potential-1} in $D^\pm$ and the Rankine-Hugoniot
condition on $S$:
\begin{equation}\label{FBConditionSelfSim-0}
\left[\rho(|\nabla\varphi|^2,\varphi)\nabla\varphi\cdot\nu\right]_S=0,
\end{equation}
where the bracket $[\cdot]$ denotes the difference of the values of
the quantity along the two sides of $S$.

Then the plane incident shock solution in the $(t, {\bf
x})$--coordinates with states $(\rho, \nabla_{\bf x}\Phi)=(\rho_0,
0,0)$ and $(\rho_1, u_1,0)$ corresponds to a continuous weak
solution $\varphi$ of (\ref{potential-eq}) in the self-similar
coordinates $(\xi,\eta)$ with the following form:
\begin{eqnarray}
&&\varphi_0(\xi,\eta)=-\frac{1}{2}(\xi^2+\eta^2) \qquad
 \hbox{for } \,\, \xi>\xi_0,
 \label{flatOrthSelfSimShock1} \\
&&\varphi_1(\xi,\eta)=-\frac{1}{2}(\xi^2+\eta^2)+ u_1(\xi-\xi_0)
\qquad
 \hbox{for } \,\, \xi<\xi_0,
 \label{flatOrthSelfSimShock2}
\end{eqnarray}
respectively, where
\begin{eqnarray}
&&u_1
=\sqrt{\frac{2(\rho_1-\rho_0)(\rho_1^{\gamma-1}-\rho_0^{\gamma-1})}
 {(\gamma-1)(\rho_1+\rho_0)}}>0,\label{state-1-velocity}\\
&&\xi_0=\rho_1
\sqrt{\frac{2(\rho_1^{\gamma-1}-\rho_0^{\gamma-1})}{(\gamma-1)(\rho_1^2-\rho_0^2)}}
=\frac{\rho_1u_1}{\rho_1-\rho_0}>0 \label{shocklocation}
\end{eqnarray}
are the velocity of state (1) and the location of the incident
shock, uniquely determined by $(\rho_0,\rho_1,\gamma)$ through
(\ref{FBConditionSelfSim-0}). Then
 $P_0=(\xi_0, \xi_0 \tan\theta_w)$ in Fig. 2, and
{\bf Problem 3.1} in the context of the potential flow equation can
be formulated as:

\medskip
{\bf Problem 7.1}\label{BVP} {(Boundary value problem)} (see Fig.
{\rm 2}). {\it Seek a solution $\varphi$ of equation
\eqref{potential-1} with \eqref{1.1.6} in the self-similar domain
$\Lambda$ with the boundary condition on
$\partial\Lambda$:
\begin{equation}\label{boundary-condition-3}
\nabla\varphi\cdot\nu|_{\partial\Lambda}=0,
\end{equation}
and the asymptotic boundary condition at infinity:
\begin{equation}\label{boundary-condition-2}
\varphi\to\bar{\varphi}:=
\begin{cases} \varphi_0 \qquad\mbox{for}\,\,\,
                         \xi>\xi_0, \eta>\xi \tan\theta_w,\\
              \varphi_1 \qquad \mbox{for}\,\,\,
                          \xi<\xi_0, \eta>0,
\end{cases}
\qquad \mbox{when $\xi^2+\eta^2\to \infty$},
\end{equation}
where {\rm (\ref{boundary-condition-2})} holds in the sense that
$
\displaystyle
\lim_{R\to\infty}\|\varphi-\overline{\varphi}\|_{C(\Lambda\setminus
B_R(0))}=0.
$
}

\medskip
For our problem, since $\varphi_1$ does not satisfy the slip
boundary condition \eqref{boundary-condition-3}, the solution must
differ from $\varphi_1$ in $\{\xi<\xi_0\}\cap\Lambda$, thus a shock
diffraction-diffraction by the wedge vertex occurs. In Chen-Feldman
\cite{ChenFeldman2a,ChenFeldman2b}, we first followed the von
Neumann criterion and the stability criterion introduced in Section
5 to establish a local existence theory of regular shock reflection
near the reflection point $P_0$ in the level of potential flow,
when the wedge angle is large and close to $\pi/2$.  In this case,
the vertical line is the incident shock $S=\{\xi=\xi_0\}$ that hits
the wedge at the point $P_0=(\xi_0, \xi_0 \tan\theta_w)$, and state
(0) and state (1) ahead of and behind $S$ are given by $\varphi_0$
and $\varphi_1$ defined in \eqref{flatOrthSelfSimShock1} and
\eqref{flatOrthSelfSimShock2}, respectively. The solutions $\varphi$
and $\varphi_1$ differ only in the domain $P_0\PtUpL \PtLwL \PtLwR$
because of shock diffraction by the wedge vertex, where the curve
$P_0\PtUpL \PtLwL$ is the reflected shock with the straight segment
$P_0\PtUpL$. State (2) behind $P_0\PtUpL$ can be computed explicitly
with the form:
\begin{equation}\label{state2}
\varphi_2(\xi,\eta)=-\frac{1}{2}(\xi^2+\eta^2)+u_2(\xi-\xi_0)+
(\eta-\xi_0\tan\theta_w)u_2\tan\theta_w,
\end{equation}
which satisfies $ \nabla\varphi\cdot \nu=0$ on $\partial\Lambda\cap
\{\xi>0\}$;  the constant velocity $u_2$ and the angle $\theta_s$
between $P_0\PtUpL$ and the $\xi$--axis are determined by
$(\theta_w,\rho_0,\rho_1,\gamma)$ from the two algebraic equations
expressing (\ref{FBConditionSelfSim-0}) and continuous  matching of
state (1) and state (2) across $P_0\PtUpL$, whose existence is
exactly guaranteed by the condition on
$(\theta_w,\rho_0,\rho_1,\gamma)$ under which regular shock
reflection-diffraction is expected to occur as in Theorem 5.1.
Moreover, $\varphi_2$ is the unique solution in the domain
$P_0\PtUpL \PtUpR$, as argued in \cite{ChangChen,Serre}. Denote
$$
P_1P_4:=\Sonic=\partial\Omega\cap\partial B_{c_2}(u_2,
u_2\tan\theta_w)
$$
the sonic arc of state $(2)$ with center $(u_2,
u_2\tan\theta_w)$ and radius $c_2$. Also we introduce the following
notation for the other parts of $\der\Om$:
\begin{equation*}
\Shock:=P_1P_2; \qquad \Wedge:=\der \Om \cap \der
\Lambda\cap\{\eta>0\}\equiv \PtLwR\PtUpR; \quad
\Gamma_{symm}:=\{\eta=0\}\cap\partial\Omega.
\end{equation*}

Then {\bf Problem 7.1} can be formulated as:

\medskip
{\bf Problem 7.2}. {\it Seek a solution $\varphi$ in $\Omega$ to
equation \eqref{potential-1} with \eqref{1.1.6}
subject to the boundary condition \eqref{boundary-condition-3} on
$\partial\Omega\cap\partial\Lambda$,
the Rankine-Hugoniot conditions on the shock $\Gamma_{shock}$:
\begin{eqnarray}
&&[\varphi]_{\Shock}=0, \label{RH-1}\\
&&[\rho(\nabla\varphi,\varphi,\rho_0)\nabla\varphi\cdot\nu]_{\Shock}=0,\label{RH-2}
\end{eqnarray}
and the Dirichlet boundary condition on the sonic arc
$\Gamma_{sonic}$:
\begin{equation}\label{sonic-dirichlet}
 (\varphi-\varphi_2)|_{\Gamma_{sonic}}=0.
\end{equation}
}

It should be noted that, in order that the solution $\varphi$ in the
domain $\Omega$ is a part of the global solution to {\bf Problem
7.1}, that is, $\varphi$ satisfies the equation in the sense of
distributions in $\Lambda$, especially across the sonic arc
$\Gamma_{sonic}$, it is requaired that
$$
\nabla(\varphi-\varphi_2)\cdot \nu|_{\Gamma_{sonic}}=0.
$$
That is, we have to match our solution with state (2), which is the
necessary condition for our solution in the domain $\Omega$ to be a
part of the global solution. To achieve this, we have to show that
our solution is at least $C^1$ with $\nabla (\varphi-\varphi_2)=0$
across $\Gamma_{sonic}$.

Then the problem can be reformulated as the following free boundary
problem:

\medskip
{\bf Problem 7.3 (Free boundary problem)}. {\it Seek a solution
$\varphi$ and a free boundary $\Shock=\{\xi=f(\eta)\}$ such that
\begin{itemize}
\item[(i)] $f\in C^{1,\alpha}$ and
\begin{equation}\label{7.11}
\Omega_+=\{\xi>f(\eta)\}\cap D{=\{\varphi<\varphi_1\}\cap D};
\end{equation}

\item[(ii)] $\varphi$ satisfies the free boundary condition \eqref{RH-2}
along $\Shock$;

\item[(iii)] $\varphi\in C^{1,\alpha}(\overline{\Omega_+})\cap C^2(\Omega_+)$
     solves \eqref{potential-eq} in $\Omega_+$, is subsonic in $\Omega_+$,
     and satisfies
\begin{eqnarray}
&&(\varphi-\varphi_2,
\nabla(\varphi-\varphi_2)\cdot\nu)|_{\Gamma_{sonic}}=0,
  \label{7.12}\\
&&\nabla \varphi\cdot\nu|_{\Wedge\cup\Gamma_{symm}}=0. \label{7.13}
\end{eqnarray}
\end{itemize}
}

The boundary condition on $\Gamma_{symm}$ implies that $f'(0)=0$ and
thus ensures the orthogonality of the free boundary with the
$\xi$-axis. Formulation \eqref{7.11} implies that the free boundary
is determined by the level set $\varphi=\varphi_1$, which is a
convenient formulation to apply useful free boundary techniques. The
free boundary condition \eqref{RH-2} along $\Shock$ is the conormal
boundary condition on $\Shock$. Condition \eqref{7.12} ensures that
the solution of the free boundary problem in $\Omega$ is a part of
the global solution as pointed out earlier. Condition \eqref{7.13}
is the slip boundary condition.

{\bf Problem 7.3} involves two types of transonic flow: one is a
continuous transition through the sonic arc $\Sonic$ as a fixed
boundary from the pseudo-supersonic region (2) to the
pseudo-subsonic region $\Omega$; the other is a jump transition
through the transonic shock as a free boundary from the supersonic
region (1) to the subsonic region $\Omega$.

\section{Global Theory for Regular Reflection-Diffraction
for Potential Flow}

In this section, we describe a global theory for regular
reflection-diffraction established in Chen-Feldman
\cite{ChenFeldman2a,ChenFeldman2b,ChenFeldman3} and Bae-Chen-Feldman
\cite{BCF}.

\subsection{Existence and stability of regular
reflection-diffraction configurations}

In Chen-Feldman \cite{ChenFeldman2a,ChenFeldman2b}, we have
developed a rigorous mathematical approach to solve {\bf Problem
7.3} and established a global theory for solutions of regular
reflection-diffraction, which converge to the unique solution of the
normal shock reflection when $\theta_w$ tends to $\pi/2$.

Introduce the polar coordinates $(r,\theta)$ with respect to the
center $(u_2, u_2\tan\theta_w)$ of the sonic arc $\Sonic$ of state
(2), that is,
\begin{equation}\label{coordPolar}
\xi-u_2=r\cos{\theta},\quad \eta-u_2\tan\theta_w=r\sin{\theta}.
\end{equation}
Then, for $\eps\in (0, c_2)$, we denote by
$$
\Omega_\eps:=\Omega\cap\{(r,\theta)\; : \;0<c_2-r<\eps\}
$$
the
$\eps$-neighborhood of the sonic arc $\PtUpL\PtUpR$ within $\Omega$;
see Fig. 4. In $\Omega_\eps$, we introduce the coordinates:
\begin{equation}\label{coordNearSonic}
x=c_2-r, \quad y=\theta-\theta_w.
\end{equation}
Then $\Omega_\eps\subset \{0<x<\eps, \; y>0\}$ and
$\PtUpL\PtUpR\subset \{x=0 \; y>0\}$.

\begin{theorem}[Chen-Feldman \cite{ChenFeldman2a,ChenFeldman2b}] There exist
$\theta_c=\theta_c(\rho_0,\rho_1,\gamma) \in (0,\pi/2)$ and
$\alpha=\alpha(\rho_0,\rho_1,\gamma)\in (0, 1/2)$ such that, when
$\theta_w\in (\theta_c,\pi/2)$, there exists a global self-similar
solution
$$
\Phi(t, {\bf x}) =t\,\varphi(\frac{\bf x}{t}) +\frac{|\bf x|^2}{2t}
\qquad\mbox{for}\,\, \frac{\bf x}{t}\in \Lambda,\, t>0
$$
with $\rho(t, {\bf x})=(\rho_0^{\gamma-1}-\Phi_t
      -\frac{1}{2}|\nabla_{\bf x}\Phi|^2)^{\frac{1}{\gamma-1}}$
      of {\bf Problem 7.1} (equivalently, {\bf Problem 7.2}) for shock
reflection-diffraction by the wedge. The solution $\varphi$
satisfies that, for $(\xi,\eta)={\bf x}/t$,
\begin{equation}\label{phi-states-0-1-2}
\begin{split}
&\varphi\in C^{0,1}(\Lambda),\\
&\varphi\in C^{\infty}(\Omega)\cap C^{1,\alpha}(\bar{\Omega}),\\
&\varphi=\left\{\begin{array}{ll}
\varphi_0 \qquad\mbox{for}\,\, \xi>\xi_0 \mbox{ and } \eta>\xi\tan\theta_w,\\
\varphi_1 \qquad\mbox{for}\,\, \xi<\xi_0
  \mbox{ and above the reflected shock} \,\, P_0\PtUpL\PtLwL,\\
\varphi_2 \qquad \mbox{in}\,\, P_0\PtUpL\PtUpR.
\end{array}\right.
\end{split}
\end{equation}
Moreover,
\begin{enumerate}
\renewcommand{\theenumi}{\roman{enumi}}
\item \label{ellipticityInOmega}
equation \eqref{potential-eq} is elliptic in $\Omega$;

\item \label{phi-GE-phi2}
$
\vphi_2\le \vphi \le \vphi_1\qquad\inn\;\;\Om$;

\item\label{BdryReg}
  the
reflected shock $P_0\PtUpL\PtLwL$ is $C^2$ at $\PtUpL$ and
$C^\infty$ except $\PtUpL$;

\item\label{C11NormEstimate}
there exists $\eps_0\in (0, c_2/2)$ such that $\vphi\in
C^{1,1}(\overline{\Om_{\eps_0}})\cap C^2(\overline{\Om_{\eps_0}}
\setminus \overline\Sonic)$; in particular, in the coordinates {\rm
(\ref{coordNearSonic})},
\begin{equation}\label{parabolicNorm}
\|\vphi-\vphi_2\|^{(par)}_{2,0,\Om_{\eps_0}} \!:=\!\sum_{0\le k+l
\le 2}\sup_{(x,y)\in\Om_{\eps_0}}\big(x^{k+\frac l2
-2}|\der_x^k\der_y^l(\vphi-\vphi_2)(x,y)|\big)<\infty;
\end{equation}

\item\label{psi-x-est}
there exists  $\delta_0>0$ so that, in the coordinates {\rm
(\ref{coordNearSonic})},
\begin{equation}\label{ell}
|\partial_x(\vphi-\vphi_2)(x,y)|\le
\frac{2-\delta_0}{\gam+1}x\qquad\;\inn \;\;\Om_{\eps_0};
\end{equation}

\item\label{FBnearSonic}
there exist $\omega>0$ and a function $y=\hat{f}(x)$ such that, in
the coordinates {\rm (\ref{coordNearSonic})},
\begin{equation}\label{OmegaXY}
\begin{split}
&\Om_{\eps_0}=\{(x,y)\,:\,x\in(0,\;\eps_0),\;\; 0< y<\hat{f}(x)\}, \\
&\Shock\cap\partial\Omega_{\eps_0}=\{(x,y):x\in(0,\;\eps_0),\;\;
y=\hat{f}(x)\},
\end{split}
\end{equation}
and
\begin{equation} \label{fhat}
\|\hat{f}\|_{C^{1,1}([0,\;\eps_0])}<\infty,\qquad
\frac{d\hat{f}}{dx}\!\ge\omega>0\!\;\;\;
\text{for}\;\;0<x<\eps_0.
\end{equation}
\end{enumerate}
Furthermore, the solution $\varphi$ is stable with respect to the
wedge angle $\theta_w$ in $W^{1,1}_{loc}$ and converges in
$W^{1,1}_{loc}$ to the unique solution of the normal reflection as
$\theta_w\to\pi/2$.
\end{theorem}

The existence of a solution $\varphi$ of {\bf Problem 7.1}
(equivalently, {\bf Problem 7.2}), satisfying
(\ref{phi-states-0-1-2}) and property (\ref{C11NormEstimate}),
follows from \cite[Main Theorem]{ChenFeldman2b}. Property
(\ref{ellipticityInOmega}) follows from  Lemma 5.2 and  Proposition
7.1 in \cite{ChenFeldman2b}. Property (\ref{phi-GE-phi2}) follows
from Proposition 7.1 and Section 9 in \cite{ChenFeldman2b} which
assert that $\varphi-\varphi_2\in\mathcal K$, where the set
$\mathcal K$ is defined by (5.15) in \cite{ChenFeldman2b}.
Property (\ref{psi-x-est})
follows from Propositions 8.1--8.2 and Section 9 in
\cite{ChenFeldman2b}. Property (\ref{FBnearSonic}) follows from
(5.7) and (5.25)--(5.27) in \cite{ChenFeldman2b} and the fact that
$\varphi-\varphi_2\in\mathcal K$.

We remark that estimate \eqref{parabolicNorm} above confirms that
our solutions satisfy the assumptions of Theorem 6.1 for the
velocity field $(U,V)=\nabla\varphi$.

\medskip
One of the main difficulties for the global existence is that the
ellipticity condition \eqref{1.1.8} for (\ref{potential-1}) with
\eqref{1.1.6} is hard to control, in comparison to our  work on
steady flow \cite{ChenFeldman1,ChenFeldman1a,ChenFeldman2,CFe3}. The
second difficulty is that the ellipticity degenerates along the
sonic arc $\Sonic$. The third difficulty is that, on $\Sonic$, the
solution in $\Omega$ has to be matched with $\varphi_2$ at least in
$C^1$, i.e., the two conditions on the fixed boundary $\Sonic$: the
Dirichlet and conormal conditions, which are generically
overdetermined for an elliptic equation since the conditions on the
other parts of boundary have been prescribed. Thus, one needs to
prove that, if $\varphi$ satisfies (\ref{potential-1}) in $\Omega$,
the Dirichlet continuity condition on the sonic arc, and the
appropriate conditions on the other parts of $\partial\Omega$
derived from Problem 7.3, then the normal derivative
$\nabla\varphi\cdot \nu$ automatically matches with
$\nabla\varphi_2\cdot \nu$ along $\Sonic$. Indeed, equation
(\ref{potential-1}), written in terms of the function
$\psi=\varphi-\varphi_2$ in the $(x,y)$--coordinates defined near
$\Sonic$ such that $\Sonic$ becomes a segment on $\{x=0\}$, has the
form:
\begin{equation}\label{degenerate-equation}
\big(2x-(\gamma+1)\psi_x\big)\psi_{xx}+\frac{1}{c_2^2}\psi_{yy}-\psi_x=0
 \qquad\,\, \mbox{in } x>0 \mbox{  and near } x=0,
\end{equation}
plus the {\it small} terms that are controlled by $\pi/2-\theta_w$
in appropriate norms. Equation \eqref{degenerate-equation} is {\it
elliptic} if $\psi_x<2x/(\gamma+1)$. Hence, it is required to obtain
the $C^{1,1}$ estimates near $\Sonic$ to ensure
$|\psi_x|<2x/(\gamma+1)$ which in turn implies both the ellipticity
of the equation in $\Omega$ and the match of normal derivatives
$\nabla\varphi\cdot \nu=\nabla\varphi_2\cdot\nu$ along $\Sonic$.
Taking into account the {\it small} terms to be added to equation
\eqref{degenerate-equation}, one needs to make the stronger estimate
$|\psi_x|\le 4x/\big(3(\gamma+1)\big)$ and assume that
$\pi/2-\theta_w$ is suitably small to control these additional
terms. Another issue is the non-variational structure and
nonlinearity of this problem which makes it hard to apply directly
the approaches of Caffarelli \cite{Ca} and Alt-Caffarelli-Friedman
\cite{AC,ACF}. Moreover, the elliptic degeneracy and geometry of the
problem makes it difficult to apply the hodograph transform approach
in Chen-Feldman \cite{ChenFeldman2} and Kinderlehrer-Nirenberg
\cite{KinderlehrerNirenberg} to fix the free boundary.

For these reasons, one of the new ingredients in our approach is to
develop further the iteration scheme in
\cite{ChenFeldman1,ChenFeldman2} to a partially modified equation.
We modified equation  (\ref{potential-1}) in $\Omega$ by a proper
Shiffmanization (i.e. a cutoff) that depends on the distance to the
sonic arc, so that the original and modified equations coincide when
$\varphi$ satisfies $|\psi_x| \le 4x/\big(3(\gamma+1)\big)$, and the
modified equation ${\mathcal N}\varphi=0$ is elliptic in $\Omega$
with elliptic degeneracy on $\PtUpL\PtUpR$. Then we solved a free
boundary problem for this modified equation: The free boundary is
the curve $\Shock$, and the free boundary conditions on $\Shock$ are
$\varphi=\varphi_1$ and the Rankine-Hugoniot condition
(\ref{FBConditionSelfSim-0}).

On each step, an {\it iteration free boundary} curve $\Sonic$ is
given, and  a solution of the modified equation ${\mathcal
N}\varphi=0$ is constructed in $\Omega$ with the boundary condition
(\ref{FBConditionSelfSim-0}) on $\Shock$, the Dirichlet condition
$\varphi=\varphi_2$ on the degenerate arc $\Sonic$, and
$\nabla\varphi\cdot \nu=0$ on $\PtLwL\PtLwR$ and $\Wedge$. Then we
proved that $\varphi$ is in fact $C^{1,1}$ up to the boundary
$\Sonic$, especially $|\nabla(\varphi-\varphi_2)|\le Cx$, by using
the nonlinear structure of elliptic degeneracy near $\Sonic$ which
is modeled by equation (\ref{degenerate-equation}) and a scaling
technique similar to Daskalopoulos-Hamilton \cite{DH} and Lin-Wang
\cite{LW}. Furthermore, we modified the {\it iteration free
boundary} curve $\Shock$ by using the Dirichlet condition
$\varphi=\varphi_1$ on $\Shock$.
 A fixed point $\varphi$ of this {\it iteration procedure} is a solution
of the free boundary problem for the modified equation. Moreover, we
proved the precise gradient estimate:
$|\psi_x|<4x/\big(3(\gamma+1)\big)$ for $\psi$, which implies that
$\varphi$ satisfies the original equation
(\ref{potential-eq}).

\smallskip
This global theory for large-angle wedges has been extended in
Chen-Feldman \cite{ChenFeldman3} to the sonic angle $\theta_s\le
\theta_c$, for which state (2) is sonic, such that, as long as
$\theta_w\in (\theta_s, \pi/2]$, the global regular
reflection-diffraction configuration exists.

\begin{theorem}[von Neumann's Sonic Conjecture (Chen-Feldman
\cite{ChenFeldman3})] The global existence result in Theorem {\rm
8.1} can be extended up to the sonic wedge-angle $\theta_s$ for any
$\gamma\ge 1$ and $u_1\le c_1$. Moreover, the solutions satisfy the
properties {\rm (i)}--{\rm (vi)} in Theorem {\rm 8.1}.
\end{theorem}

The condition $u_1\le c_1$ depends explicitly only on the parameters
$\gamma>1$ and $\rho_1>\rho_0>0$.  For the case $u_1>c_1$, we have
been making substantial progress as well, and the final detailed
results can be found in Chen-Feldman \cite{ChenFeldman3}.

\subsection{Optimal regularity}

By Theorem 8.1(\ref{C11NormEstimate}), the solution $\varphi$
constructed there is at least $C^{1,1}$ near the sonic arc $\Sonic$.
The next question is to analyze the behavior of solutions
$\varphi(\xi,\eta)$ to regular reflection-diffraction, especially
the optimal regularity of the solutions.

\smallskip
We first define the class of regular
reflection-diffraction solutions.

\begin{definition}\label{RegReflSolDef}
{\it Let $\gamma > 1$, $\rho_1>\rho_0>0$, and $\theta_w\in(0,
\pi/2)$ be constants, let  $u_1$ and $\xi_0$ be defined by {\rm
(\ref{state-1-velocity})} and {\rm (\ref{shocklocation})}.
Let the incident shock $S=\{\xi=\xi_0\}$ hits the wedge at the point
$P_0=(\xi_0, \xi_0 \tan\theta_w)$, and let state $(0)$ and state
$(1)$ ahead of and behind $\Shock$ be given by
\eqref{flatOrthSelfSimShock1} and \eqref{flatOrthSelfSimShock2},
respectively.
The function $\vphi\in C^{0,1}(\Lambda)$ is a regular
reflection-diffraction solution if $\vphi$ is a solution to {\bf
Problem 7.1}
such that

{\rm (a)} there exists state $(2)$ of form {\rm (\ref{state2})} with
$u_2>0$, satisfying the entropy condition $\rho_2>\rho_1$ and the
Rankine-Hugoniot condition  {\rm (\ref{FBConditionSelfSim-0})} along
the line $S_1:=\{\varphi_1=\varphi_2\}$ which contains the points
$P_0$ and $\PtUpL$, such that $\PtUpL\in \Lambda$ is on the sonic
circle of state $(2)$, and state $(2)$ is supersonic along
$P_0\PtUpL$;

{\rm (b)} there exists an open, connected domain
$\Omega:=P_1P_2P_3P_4\subset \Lambda$ such that {\rm
(\ref{phi-states-0-1-2})} holds and  equation \eqref{potential-eq}
is elliptic in $\Omega$;

{\rm (c)} $\varphi\ge \varphi_2$ on the part $P_1P_2=\Gamma_{shock}$
of the reflected shock.}
\end{definition}

\begin{remark}
The global solution constructed in
\cite{ChenFeldman2a,ChenFeldman2b,ChenFeldman3} is a regular
reflection-diffraction solution, which is a part of the assertions
in Theorems 8.1--8.2.
\end{remark}

\begin{remark}
If state $(2)$ exists and is supersonic, then the line
$S_1=\{\varphi_1=\varphi_2\}$ necessarily intersects the sonic
circle of state $(2)$; see the argument in
\cite{ChenFeldman2a,ChenFeldman2b} starting from {\rm (3.5)} there.
Thus, the only assumption regarding the point $\PtUpL$ is that $S_1$
intersects the sonic circle within $\Lambda$.
\end{remark}

\begin{remark}
We note that, in the case $\theta_w=\frac{\pi}2$, the regular
reflection becomes the normal reflection, in which $u_2=0$ and the
solution is smooth across the sonic line of state {\rm (2)}; see
\cite[Section 3.1]{ChenFeldman2b}. Condition $\theta_w\in(0,
\frac{\pi}2)$ in Definition {\rm \ref{RegReflSolDef}} rules out this
case. Moreover, for $\theta_w\in(0, \frac{\pi}2)$,
 the property $u_2>0$ in part {\rm (a)} of Definition {\rm
 \ref{RegReflSolDef}}
is always true for state $(2)$ of form {\rm (\ref{state2})},
satisfying the entropy condition $\rho_2>\rho_1$ and the
Rankine-Hugoniot condition (7.1) along the line
$S_1:=\{\varphi_1=\varphi_2\}$ which contains the point $P_0$. These
are readily derived from the calculations in \cite[Section
3.2]{ChenFeldman2b}.
\end{remark}

\begin{remark}
There may exist a global regular reflection-diffraction
configuration when state $(2)$ is subsonic which is a very narrow
regime \cite{CF,Neumann1,Neumann2}. Such a case does not involve the
difficulty of elliptic degeneracy, which we are facing for the
configurations in the class of solutions in the sense of Definition
8.3.
\end{remark}

\begin{remark}
Since $\varphi=\varphi_1$ on $\Shock$ by \eqref{phi-states-0-1-2},
condition {\rm (c)} in Definition {\rm \ref{RegReflSolDef}} is
equivalent to
$$
\Gamma_{shock}\subset \{\varphi_2\le \varphi_1\},
$$
that is, $\Gamma_{shock}$ is below $S_1$.
\end{remark}

In Bae-Chen-Feldman \cite{BCF}, we have developed a mathematical
approach
to establish the regularity of solutions of the regular
reflection-diffraction problem in the sense of Definition 8.3.

First, we have shown that any regular reflection-diffraction
solutions cannot be $C^2$ across the sonic arc  $\Sonic:=\PtUpL
\PtUpR$.

\begin{theorem}[Bae-Chen-Feldman
\cite{BCF}] \label{noC2acrossSonicLineThm} There does not exist a
global regular reflection-diffraction solution in the sense of
Definition {\rm 8.3} such that $\vphi$ is $C^2$ across the sonic arc
$\Sonic$.
\end{theorem}

Now we study the one-sided regularity up to $\Sonic$ from the
elliptic side, i.e., from $\Omega$. For simplicity of presentation,
we now use a localized version of $\Omega_\varepsilon$: For a given
neighborhood $\mathcal{N}(\Gamma_{sonic})$ of $\Gamma_{sonic}$ and
$\varepsilon>0$, define
$$
\Omega_\varepsilon :=\Omega\cap
\mathcal{N}(\Gamma_{sonic})\cap\{x<\varepsilon\}.
$$
Since $\mathcal{N}(\Gamma_{sonic})$ will be fixed in the following
theorem, we do not specify the dependence of $\Omega_\varepsilon$ on
$\mathcal{N}(\Gamma_{sonic})$.

\begin{theorem}[Bae-Chen-Feldman
\cite{BCF}] \label{mainthm1}
Let
$\vphi$ be a regular reflection-diffraction solution in the sense of
Definition
{\rm \ref{RegReflSolDef}} and satisfy the following properties:
There exists a neighborhood $\mathcal{N}(\Gamma_{sonic})$ of
$\Gamma_{sonic}$ such that
\begin{enumerate}
\item[(a)]\label{C11NormEstimate-a}
$\varphi$ is $C^{1,1}$ across the sonic arc $\Sonic$:
$\vphi\in C^{1,1}(\overline{P_0\PtUpL\PtLwL\PtLwR}\cap
{\mathcal{N}(\Sonic)})$;

\item[(b)]\label{psi-x-est-a}
there exists  $\delta_0>0$ so that, in the coordinates {\rm
(\ref{coordNearSonic})},
\begin{equation}\label{ell-a}
|\partial_x(\vphi-\vphi_2)(x,y)|\le
\frac{2-\delta_0}{\gam+1}x\;\qquad\inn \;\;\Om \cap
{\mathcal{N}(\Sonic)};
\end{equation}

\item[(c)]\label{FBnearSonic-a}
there exists $\eps_0>0$, $\omega>0$, and a function $y=\hat{f}(x)$
such that, in the coordinates {\rm (\ref{coordNearSonic})},
\begin{equation}\label{OmegaXY-a}
\begin{split}
&\Om_{\eps_0}=\{(x,y)\,:\, x\in(0,\;\eps_0),\;\; 0< y<\hat{f}(x)\}, \\
&\Shock\cap\partial\Omega_{\eps_0}=\{(x,y)\,:\,x\in(0,\;\eps_0),\;\;
y=\hat{f}(x)\},
\end{split}
\end{equation}
and
\begin{equation} \label{fhat-a}
\|\hat{f}\|_{C^{1,1}([0,\;\eps_0])}<\infty,\;\;\;\;
\frac{d\hat{f}}{dx}\!\ge\omega>0\!\;\;\; \text{for}\;\;0<x<\eps_0.
\end{equation}
\end{enumerate}
Then we have
\begin{enumerate}
\renewcommand{\theenumi}{\roman{enumi}}
\item\label{C2alpSonic}
$\vphi$ is $C^{2,\alp}$ up to $\Sonic$ away from $P_1$ for any
$\alp\in (0,1)$.
That is, for any $\alpha\in (0,1)$ and any given
$(\xi_0,\eta_0)\in\overline\Sonic\setminus\{\PtUpL\}$, there exists
$K<\infty$ depending only on $\rho_0,\,\rho_1,\,\gam,\,\eps_0,\,
\alpha, \|\vphi\|_{C^{1,1}(\Om_{\eps_0})}$, and
$d=dist((\xi_0,\eta_0),\;\Shock)$ so that
\begin{equation*}
\|\vphi\|_{2,\alp;\overline{B_{d/2}(\xi_0,\eta_0)\cap \Om}}\le K;
\end{equation*}

\item\label{limitSonic}
For any $(\xi_0,\eta_0)\in \Sonic\setminus\{\PtUpL\}$,
\begin{equation*}
\lim_{(\xi,\eta)\to (\xi_0,\eta_0)\atop
 (\xi,\eta)\in\Om }(D_{rr}\vphi-D_{rr}\vphi_2)=\frac{1}{\gamma+1};
\end{equation*}

\item\label{2ndderJumpSonic}
$D^2\vphi$ has a jump across $\Sonic$: For any $(\xi_0,\eta_0)\in
\Sonic\setminus\{\PtUpL\}$,
\begin{eqnarray*}
&&\lim_{(\xi,\eta)\to (\xi_0,\eta_0)\atop
 (\xi,\eta)\in\Om }D_{rr}\vphi \;-\;
 \lim_{(\xi,\eta)\to (\xi_0,\eta_0)\atop
 (\xi,\eta)\in\Lambda\setminus\Om }D_{rr}\vphi\;=\;
 \frac{1}{\gamma+1},\\
 &&\lim_{(\xi,\eta)\to (\xi_0,\eta_0)\atop
 (\xi,\eta)\in\Om }(D_{r\theta}, D_{\theta\theta})\vphi=
 \lim_{(\xi,\eta)\to (\xi_0,\eta_0)\atop
 (\xi,\eta)\in\Lambda\setminus\Om }(D_{r\theta}, D_{\theta\theta})\vphi\;=0;
\end{eqnarray*}

\item\label{noLimAtP1}
The limit $\lim_{(\xi,\eta)\to\PtUpL \atop (\xi,\eta)\in \Om}
D^2\vphi$ does not exist.
\end{enumerate}
\end{theorem}

We remark that the solutions established in
\cite{ChenFeldman2a,ChenFeldman2b,ChenFeldman3} satisfy the
assumptions of Theorem 8.5.
In particular, we proved that the
$C^{1,1}$-regularity is {\it optimal} for the solution across the
open part $P_1P_4$ of the sonic arc (the degenerate elliptic curve)
and at the point $P_1$ where the sonic circle meets the reflected
shock (as a free boundary).

To achieve the optimal regularity, one of the main difficulties is
that the sonic arc $\Sonic$ is the transonic boundary separating the
elliptic region from the hyperbolic region, near where the solution
is governed by the nonlinear degenerate elliptic equation
\eqref{degenerate-equation} for $\psi=\varphi-\varphi_2$.
We carefully analyzed the features of equation
\eqref{degenerate-equation} and established the $C^{2,\alpha}$
regularity of solutions in the elliptic region up to the open sonic
arc $P_1P_4$. As a corollary, we showed that the
$C^{1,1}$-regularity is actually optimal across the transonic
boundary $P_1P_2$ from the elliptic to hyperbolic region. Since the
reflected shock $P_1P_2$ is regarded as a free boundary connecting
the hyperbolic region (1) with the elliptic region $\Omega$ for the
nonlinear second-order equation of mixed type, another difficulty
for the optimal regularity of the solution is that the point $P_1$
is exactly the point where the degenerate elliptic arc $P_1P_4$
meets a transonic free boundary for the nonlinear partial
differential equation of second order. As far as we know, this is
the first optimal regularity result for solutions to a free boundary
problem of nonlinear degenerate elliptic equations at the point
where an elliptic degenerate curve meets the free boundary. To
achieve this, we carefully constructed two sequences of points on
where the corresponding sequences of values of $\psi_{xx}$ have
different limits at $P_1$; this has been done by employing the
one-sided $C^{2,\alpha}$ regularity of the solution up to
the open arc $P_1P_4$ and by studying detailed features of the free
boundary conditions on the free boundary $P_1P_2$, i.e., the
Rankine-Hugoniot conditions.

\smallskip
We remark that some efforts were also made mathematically for the
reflection-diffraction problem via simplified models. One of these
models, the unsteady transonic small-disturbance (UTSD) equation,
was derived and used in Keller-Blank \cite{KB}, Hunter-Keller
\cite{HK}, Hunter \cite{hunter1}, and Morawetz \cite{Morawetz2} for
asymptotic analysis of shock reflection-diffraction. Also see Zheng
\cite{Zheng1} for the pressure gradient equation and
Canic-Keyfitz-Kim \cite{CKK1} for the UTSD equation and the
nonlinear wave system. Furthermore, in order to deal with the
reflection-diffraction problem, some asymptotic methods have been
also developed. Lighthill \cite{Lighthill} studied shock
reflection-diffraction under the assumption that the wedge angle is
either very small or close to $\pi/2$. Keller-Blank \cite{KB},
Hunter-Keller \cite{HK}, Harabetian \cite{Harabetian}, and
Gamba-Rosales-Tabak \cite{GRT} considered the problem under the
assumption that the shock is so weak that its motion can be
approximated by an acoustic wave. For a weak incident shock and a
wedge with small angle for potential flow, by taking the jump of the
incident shock as a small parameter, the nature of the shock
reflection-diffraction pattern was explored in Morawetz
\cite{Morawetz2} by a number of different scalings, a study of mixed
equations, and matching the asymptotics for the different scalings.
Also see Chen \cite{Sxchen} for a linear approximation of shock
reflection-diffraction when the wedge angle is close to $\pi/2$
 and Serre
\cite{Serre} for an apriori analysis of solutions of shock
reflection-diffraction and related discussions in the context of the
Euler equations for isentropic and adiabatic fluids.

\medskip
Another related recent effort has been on various important physical
problems in steady potential flow, as well as steady fully Euler
flow, for which great progress has been made. The problems for
global subsonic flow past an obstacle and for local supersonic flow
past an obstacle with sharp head are classical, due to the works of
Shiffman \cite{Sh}, Bers \cite{Bers}, Finn-Gilbarg
\cite{FinnGilbarg}, Dong \cite{Dong}, and others for the pure
elliptic case and to the works of Gu \cite{Gu}, Shaeffer
\cite{Schaeffer}, Li \cite{Lit},
and others for the pure hyperbolic case. Recent progress has been
made on transonic flow past nozzles (e.g.
\cite{BaeFeldman,CCSo,CCF,ChenFeldman1,ChenFeldman1a,CFe3,CY,Sxchen3,CSY,XY,Yuan1,Yuan2}),
transonic flow past a wedge or conical body (e.g. \cite{C-Fang,
ChenFang,Fang}), and transonic flow past a smooth obstacle (e.g.
\cite{CSW1,GM1,Mora1,Mora3}). Also see \cite{CDSD,XX} for the
existence of global subsonic-sonic flow,
\cite{CZZ:1,ChS1,ChS2,ChS3,ChS4,CXY,Lie,Zh1} for global supersonic
flow past an obstacle with sharp head, and the references cited
therein.

For some of other recent related developments, we refer the reader
to Chen \cite{SxchenM-1,SxchenM-2} for the local stability of Mach
configuration, Elling-Liu \cite{EllingLiu2} for physicality of weak
Prandtl-Meyer reflection for supersonic potential flow around a
ramp, Serre \cite{Serre08} for multidimensional shock interaction
for a Chaplygin gas, Canic-Keyfitz-Kim \cite{CKK1,CKK2} for
semi-global solutions for the shock reflection problem,
Glimm-Ji-Li-Li-Zhang-Zhang-Zheng \cite{GJLZ} for the formation of a
transonic shock in a rarefaction Riemann problem for polytropic
gases, Zheng \cite{Zhe,Zhe1,Zheng1} for various solutions to some
two-dimensional Riemann problems, Gues-M\'{e}tivier-Williams-Zumbrun
\cite{GMWZ1,GMWZ2} and Benzoni-Gavage and Serre \cite{BG} for
revisits of the local stability of multidimensional shocks and phase
boundaries, among many others.

\section{Shock Reflection-Diffraction vs New Mathematics}

As we have seen from the previous discussion, the shock
reflection-diffraction problem involves several core challenging
difficulties that we have to face in the study of nonlinear partial
differential equations. These nonlinear difficulties include {\it
free boundary problems, oblique derivative problems for nonsmooth
domains, degenerate elliptic equations, degenerate hyperbolic
equations, transport equations with rough coefficients, mixed and/or
composite equations of hyperbolic-elliptic type, behavior of
solutions when a free boundary meets an elliptic degenerate curve,
and compressible vortex sheets}.
More efficient numerical methods are also required for further
understanding of shock reflection-diffraction phenomena.

Furthermore, the wave patterns of shock reflection-diffraction
configurations are the core patterns and configurations for the
global solutions of the two-dimensional Riemann problem; these
solutions are building blocks and local structure of general entropy
solutions and determine global attractors and asymptotic states of
entropy solutions, as time goes infinity, for two-dimensional
systems of  hyperbolic conservation laws.

Therefore, a successful solution to the shock reflection-diffraction
problem not only provides our understanding of shock
reflection-diffraction phenomena and behavior of entropy solutions
to multidimensional conservation laws, but also provides important
new ideas, insights, techniques, and approaches for our developments
of more efficient analytical techniques and methods to overcome the
core challenging difficulties in multidimensional problems in
conservation laws and other areas in nonlinear partial differential
equations. The shock reflection-diffraction problem is also an
excellent test problem to examine our capacity and ability to solve
rigorously various challenging problems for nonlinear partial
differential equations and related applications.

\bibliographystyle{amsalpha}

\end{document}